# Numerical-Analytical Investigation into Impact Pipe Driving in Soil with Dry Friction.
# Part II: Deformable External Medium

N. I. Aleksandrova

*N.A. Chinakal Institute of Mining, Siberian Branch, Russian Academy of Sciences,
Krasnyi pr. 54, Novosibirsk, 630091 Russia
e-mail: alex@math.nsc.ru*

**Abstract**—Under analysis is travel of P-waves in an elastic pipe partly embedded in soil with dry friction. The mathematical formulation of the problem on impact pipe driving in soil is based on the model of axial vibration of an elastic bar, considering lateral resistance described using the law of solid dry friction. The author solves problems on axial load on pipe in interaction with external elastic medium, and compares the analytical and numerical results obtained with and without accounting for the external medium deformability.

The problems on the behavior of pipes in soil appear in many engineering procedures such as driving and pulling of piles, trenchless laying of underground utilities, underground pipework laying. Motion of various rods in mechanical systems always goes with friction. It is important to understand the influence exerted on wave process by friction of external environment and side surface of a pipe or a rod. Research of the interaction of bodies with friction is extensive [1–40].

The wave processes in rods embedded in a medium were first considered by Gersevanov [11] in 1930, based on Saint-Venant's theory. Later on, propagation of P-waves in a rod with constant-amplitude and different-oriented dry friction was an object of study [7–25]. These studies assumed the rod was driven in undeformable medium. A general approach to solving such problems is reducing the nonlinear problem to a set of successively solved linear problems. It is required to define boundaries in time and axial coordinate for domains where friction has constant value and orientation. The axis of friction depends on the sign of the wave velocity. Solution in each of the domains is constructed considering the preliminary determined sign of the velocity. This approach allowed solving a lot of problems with sufficiently simple laws of loading exerted on the rod.

Nonstationary and stationary one-dimensional problems on pipe and soil interaction by the law of external dry friction in the framework of the model of sonic and ultrasonic flow of a finite and infinite rod can be found in [8, 13, 15], and two-dimensional problem of pipe driving in soil is discussed in [19].

Many researches focused on problems on piling [26–35]. In [26] experimental data on hammer accelerations are used to calculate distribution of soil resistance along a pile. In [27] a simple dynamic model is offered to find soil resistance coefficient in the course of piling; the calculated results are compared with the dynamic testing database, and it is discussed how to use the soil resistance model in the wave equation of the pile motion. Dynamic analysis of side surface of pile uses boundary element approach in [28], and the finite element method is chosen to solve axially symmetric two-dimensional problem on piling, with algorithm of slip friction contact on the pile and soil interface in [29, 30]. Experiments on influence exerted by soil plug inside an open-end driven pile on the static and dynamic response of the pile are described in [31]. In [32], based on analysis of the wave equation from [22], the load-bearing capacity of piles is calculated and compared with the experimental results. The problem on piling in [33, 34] is solved with modeling the friction contact between the pile and soil using the plastic strengthening/weakening theory. The finite element analysis in [35] uses the Coulomb law of the soil and pile contact friction, and the estimates are compared with the experiment data. Many experimental approaches to pipe driving in soil are reviewed in [36–40]. However, the researchers [26–40] left the nonstationary problems on the pile and soil interaction aside of their attention.

The dynamic problem on wave propagation in a composite material composed of fiber and binder with dry friction pre-set ate their interface is discussed in [41]. Motion of the fiber and binder is described with one-dimensional wave equations solved analytically and using finite element method. The load was assumed linearly growing in time. The main attention was paid to the perturbance delay as against the wave propagation.

The given study continues the research in [21], considering deformation of the surrounding medium, compares the results of the both models and identifies the applicability of the simple model. The author uses the end-to-end numerical calculation of wave processes is systems with dry friction, that allows analyzing any type load. The obtained approximate analytical solutions on a one-dimensional problem on radially oriented perturbance propagation on slip-free elastic medium and pipe interaction and axial load applied to the pipe end are the test patterns for the finite difference solutions of two-dimensional problems with pre-set elastic interaction or dry friction on the pipe and medium contact.

## 1. PROBLEM FORMULATION

The mathematical formulation of the problem on pipe driving in soil is based on the model of longitudinal oscillations of elastic rod, considering lateral resistance. We solve the mixed dynamic axially symmetric problem on interaction of an elastic cylinder pipe and an elastic soil layer. The pipe end is applied with a longitudinal impact. The pipe motion is described with a one-dimensional wave equation, lateral resistance of soil with low strains—using simplified models of elastic medium, and yielding—by the dry friction law.

An elastic tubular rod ($R$—radius; $h$—wall thickness) with length $L$ is embedde in soil for the length $L_1$. The applied longitudinal half-sine impact impulse has amplitude $P_0$ and duration $t_0$:

$$Q(t) = P_0 \sin(\omega_* t) H_0(t_0 - t) H_0(t), \quad \omega_* = \pi/t_0, \qquad (1)$$

where $H_0(t)$ is Heaviside's step function.

Let the origin of coordinate system be the impacted end of the rod, and the $Z$ axis be oriented in parallel to the driven rod axis. The rod motion is described with the one-dimensional wave equation in terms of мещений $U(t,z)$:

$$\ddot{U} = c^2 U''_{,zz} + \tau(z, \varepsilon_{rz}). \qquad (2)$$

Here, $c = \sqrt{E/\rho}$ is P-wave velocity in the rod; $E$ is Young's modulus; $\rho$ is the rod material density; $\tau$ is response of the medium.

Initial conditions are zero:

$$U\big|_{t=0} = 0, \quad \dot{U}\big|_{t=0} = 0. \qquad (3)$$

At the rod end $z = 0$ the stress is assigned, the rod end $z = L$ is stress-free:

$$ES_t U'_{,z}\big|_{z=0} = -Q(t), \quad ES_t U'_{,z}\big|_{z=L} = 0, \qquad (4)$$

where $S_t = \pi h(2R - h)$ is the tubular rod cross-section area.

NUMERICAL-ANALYTICAL INVESTIGATION INTO IMPACT PIPE DRIVING IN SOIL

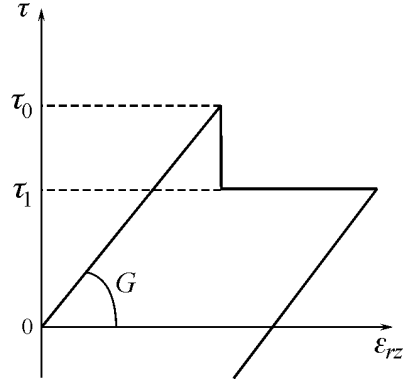

**Fig. 1.** Relation of response of the medium and shear strain.

The curve $\tau(z,\varepsilon_{rz})$ fits in the elastoplastic strain diagram, taking account of yielding peak (Fig. 1). The diagram contains an elastic line from 0 to $\tau_0$ showing the rod and medium adhesion, then ultimate shear stress $\tau_0$ changes to $\tau_1$ and frictional slip start (usually proportional to pressing force). When the current rod cross-section $z > L - L_1$ stops, the rod and medium adhere again, and the adherence can later be broken by stress waves arrived from the rod ends. Orientation of friction forces is governed by the sign of velocity of these waves. The present study assumes $\tau_0 = \tau_1$.

At the elastic stage when $\mathrm{abs}(\tau) < \tau_0$, the shear stress of soil is found using a two-dimensional model of the medium with single displacement $V$ in the line of $z$:

$$\ddot{V} = a^2 V''_{,zz} + b^2\left(V''_{,rr} + \frac{1}{r}V'_{,r}\right), \quad a^2 = \frac{\lambda + 2G}{\gamma}, \quad b^2 = \frac{G}{\gamma}. \tag{5}$$

Here, $G$ is shear modulus; $\lambda$ is Lamé's coefficient; $\gamma$ is soil density; $r$ is radial coordinate. Equation (5) is a model of the medium, including compressibility of soil along the rod axis (coefficient $a^2$), shearing resistance (coefficient $b^2$) and inertance of soil particles moving along the rod axis.

This model was used to solve many dynamic and static problems in [42], and it was shown that amplitudes of axial strains and their partial derivatives are much high than these values in radial direction. This allows transfer from the accurate model of the elastic theory to the physically equal and mathematically simpler model of medium deformation with the single displacement (5).

Boundary condition for (5) in the absence of friction slip on the rod surface is:

$$V(z,r)\big|_{r=R} = U(z), \quad L - L_1 \leq z \leq L. \tag{6}$$

It defines the equality of the displacements of the rod and soil on the rod surface. The other boundary conditions are the absence of displacements at the exterior boundary $r = R_2$ of the soil and the absence of stresses on the free surfaces of the soil:

$$V(z,r)\big|_{r=R_2} = 0 \quad (L - L_1 \leq z \leq L), \quad V'_{,z}(z,r)\big|_{z=L-L_1} = 0, \quad V'_{,z}(z,r)\big|_{z=L} = 0. \tag{7}$$

The soil response $\tau$ is:

$$\tau(z,t) = \frac{P_t}{S_t \rho}\begin{cases} GV'_{,r}\big|_{r=R}, & z > L - L_1, \\ 0, & z \leq L - L_1. \end{cases} \tag{8}$$

Here, $P_t = 2\pi R$ is perimeter of the tubular rod.

At the stage of slip when $\mathrm{abs}(\tau) > \tau_0$, the rod and soil interaction is described by the dry friction law:

$$\tau(z,t) = \frac{P_t}{S_t \rho} \begin{cases} -k\tau_0, & z > L-L_1, \\ 0, & z \leq L-L_1, \end{cases} \qquad (9)$$

$$k = \mathrm{sign}[\dot{U}(z) - \dot{V}(z,R)] \quad \text{at} \quad \dot{U} - \dot{V} \neq 0.$$

During slip, friction forces on the rod and the medium agree in value and differ in orientation, therefore, there are the following boundary conditions for the medium:

$$GV'_{,r}\big|_{r=R} = k\tau_0. \qquad (10)$$

The method of end-to-end calculation of wave processes in systems with dry friction, presented below, allows considering any type load as the friction value and orientation at any point of the rod and at any time moment is chosen based on physical considerations. Also, the author presents analytical estimates for some problems.

## 2. ANALYTICAL ALGORITHMS

The system of equations 91)–(10) was solved in the case of a unit impact, using cross-type finite difference scheme and method of minimization of data scatter. Equation (5) is approximated as:

$$V_{ji}^{n+1} - 2V_{ji}^n + V_{ji}^{n-1} = \left[a^2 \Lambda_{jj} V_{ji}^n + b^2 \left(\Lambda_{ii} V_{ji}^n + \frac{1}{r}\Lambda_i^0 V_{ji}^n\right)\right] h_t^2, \qquad (11)$$

where $\Lambda_{jj}$, $\Lambda_{ii}$ are second-order central-difference operators in coordinates $z$, $r$; $\Lambda_i^0$ is central-difference operator of the first derivative with respect to $r$; $V_{ji}^n = V(h_t n, h_z j, r_i)$ are mesh values of displacements at time $t = h_t n$ at the point $z = h_z j$; $r_i = R + h_r(i-1)$, $i = 1,2,...$; $h_t$, $h_z$, $h_r$ are mesh steps in coordinates $t$, $z$, $r$.

The approximation of (2) has the form:

$$U_j^{n+1} - 2U_j^n + U_j^{n-1} = (c^2 \Lambda_{jj} U_j^n + \tau_j^n) h_t^2, \qquad (12)$$

where $U_j^n = U(h_t n, h_z j)$, $\tau_j^n = \tau(h_t n, h_z j)$ are mesh values of displacements and shear stresses on the rod surface.

Conditions of stability of Eqs. (11) and (12) are:

$$h_t \leq \left(\frac{a^2}{h_z^2} + \frac{b^2}{h_r^2}\right)^{-1/2}, \quad h_t \leq \frac{h_z}{c}.$$

The optimized parameters of the finite difference mesh to make data scatter of (11) and (12) minimum are $h_z = ch_t$ axially and $h_r = h_t cb/\sqrt{c^2 - a^2}$ radially.

Approximation of boundary conditions used right-hand side for two derivatives at the rod end $z = 0$: $ES_t(U_1^n - U_0^n)/h_z = -Q(t)$ and central-difference approximation at $z = L$. Similarly, we approximated boundary conditions for soil at the boundaries $z = L - L_1$ and $z = L$. For the boundary condition (6), three-point approximation was used [43].

The algorithm with friction is analogous to that in [19, 21]. Neither orientation nor value of friction force are known beforehand, and we calculate relative velocities of displacing points of the medium and rod for two optional signs of $k$ ($k > 0$ and $k < 0$):

(a) When $k > 0$, introduce a fictitious velocity $\dot{U}_j^+ - \dot{V}_{j,1}^+$, where:

$$\dot{V}_{j,1}^+ = \frac{V_{j,1}^{n+1} - V_{j,1}^n + h_t^2 \tau_0}{h_t}, \quad \dot{U}_j^+ = \frac{U_j^{n+1} - U_j^n + h_t^2 \tau_0}{h_t};$$

(b) When $k < 0$, introduce a fictitious velocity $\dot{U}_j^- - \dot{V}_{j,1}^-$, where:

$$V_{j,1}^- = \frac{V_{j,1}^{n+1} - V_{j,1}^n - h_t^2 \tau_0}{h_t}, \quad U_j^- = \frac{U_j^{n+1} - U_j^n - h_t^2 \tau_0}{h_t}.$$

The values of $V_{j,1}^{n+1}$ and $U_j^{n+1}$ are calculated from (11) and (12) without regard to friction, and there are two possible situations:

1. If $\dot{U}_j^+ - \dot{V}_{j,1}^+$ and $\dot{U}_j^- - \dot{V}_{j,1}^-$ have the same sign, the true displacements $V_{j,1}^{n+1}$, $U_j^{n+1}$ out of $V_{j,1}^+$, $U_j^+$ and $V_{j,1}^-$, $U_j^-$ are selected as $V_{j,1}^k$, $U_j^k$ such that:

$$\mathrm{abs}(\dot{U}_j^k - \dot{V}_{j,1}^k) = \min[\mathrm{abs}(\dot{U}_j^- - \dot{V}_{j,1}^-), \mathrm{abs}(\dot{U}_j^+ - \dot{V}_{j,1}^+)].$$

2. If $\dot{U}_j^+ - \dot{V}_{j,1}^+$ and $\dot{U}_j^- - \dot{V}_{j,1}^-$ have different sign, or one of these velocities is zero, then, assuming friction is passive, the real relative velocity of the soil and rod equals zero. The soil and rod adhere and are moved together, thus, friction is absent, and the rod and soil interaction is elastic.

Accordingly, the problem on finding the interface of the relative motion and relative rest domains, which is the major difficulty in analytical solutions, reduces to finding points where $\dot{V}_{j,1}^+ - \dot{U}_j^+$ and $\dot{V}_{j,1}^- - \dot{U}_j^-$ have different signs, or where one of them is zero.

The value and orientation of friction force are unambiguously defined in the calculation process, and, therefore, on each time level, we solve a linear problem with defined friction included as load in the right-hand side of equation.

Friction can be either passive or active. For instance, when particles of the soil and rod are moved in the same direction and $|\dot{V}| > |\dot{U}|$, friction accelerates the rod; otherwise, it decelerates the rod. With the same sign $\dot{V}$ and $\dot{U}$, it is deceleration that takes place.

### 3. ANALYTICAL ESTIMATES

The one-dimensional problem includes a thin soil layer (in the line of $z$) with radius $R_2$ and a pipe with radius $R$ driven in soil. Motion of the medium is described with the one-dimensional wave equation in terms of the radial coordinate $r$ driven from (5) under assumption on neglectability of the derivatives with respect to the axis $z$:

$$\ddot{V} = b^2 \left( V''_{,rr} + \frac{1}{r} V'_{,r} \right), \quad b^2 = \frac{G}{\gamma}. \qquad (13)$$

The boundary conditions:

$$V(r)|_{r=R} = U, \quad V(r)|_{r=R_2} = 0. \qquad (14)$$

The pipe motion equation is derived from (2) on the assumption that the relation of the displacement and $z$ can be disregarded:

$$\ddot{U} = \frac{P_t G}{S_t \rho} V'_{,r} \bigg|_{r=R} + \frac{Q(t)}{S_t \rho}. \qquad (15)$$

The initial conditions are zero.

Take the Laplace transform of time with parameter $p$:

$$p^2 V^L = b^2 \left[ (V''_{,rr})^L + \frac{1}{r}(V'_{,r})^L \right], \quad p^2 V^L \big|_{r=R} = \frac{P_t G}{S_t \rho}(V'_{,r})^L \bigg|_{r=R} + \frac{Q^L}{S_t \rho}, \quad V^L \big|_{r=R_2} = 0. \qquad (16)$$

The first case to be analyzed is the infinite medium ($R_2 \to \infty$). Let the pipe be applied with a step load $Q(t) = P_0 H_0(t)$. Using Laplace transform of a function re-writes (16) as:

$$U^L = \frac{P_0 K_0(\eta)}{LS_t \rho p^2 (pK_0(\eta) + K_1(\eta) P_t \gamma b / S_t \rho)},$$

where $K_0$, $K_1$ are cylindrical functions of imaginary argument; $\eta = Rp/b$.

Asymptotic form of the Laplace transform of function when $p \to 0$, which is equal to $t \to \infty$ in the space of pro-images, is:

$$U^L = \frac{P_0}{2\pi Gp} \ln\left(\frac{pCR}{2b}\right).$$

Converse of the Laplace transform [44] produces the asymptotic relation of the pipe displacement and time when $t \to \infty$ in the infinite medium:

$$U(t) = \frac{P_0}{2\pi G} \ln\left(\frac{2bt}{R}\right), \quad (t \to \infty). \tag{17}$$

Considering wave reflection from the external boundary $r = R_2$, image solution of (16) is:

$$U^L = \frac{P_0}{S_t \rho p^2}\left(p + \frac{[K_0(\eta)I_0(\eta_2) - I_0(\eta)K_0(\eta_2)]}{[K_1(\eta)I_0(\eta_2) - I_1(\eta)K_0(\eta_2)]} \frac{P_t \gamma b}{S_t \rho}\right)^{-1}. \tag{18}$$

Here, $I_0$, $I_1$ are modified Bessel functions; $\eta = Rp/b$; $\eta_2 = R_2 p/b$.

If $t \to \infty$ ($p \to 0$), then it follows from (18) that the solution oscillates about relatively static equilibrium defined by the formula:

$$U_{stat} = \frac{P_0}{2\pi G} \ln\left(\frac{R_2}{R}\right). \tag{19}$$

Figure 2 shows the results of the finite difference solution of the system (13)–(15) by solid lines, asymptotic solutions of (17) by dashed lines and the static value of (19) by dash-and-dotted line. The thick lines correspond to $R_2 = 20$ m, the thin lines— $R_2 = 2$ m. the finite difference mesh parameters were: $h_r = 0.01$ m, $h_t = h_r/b$. The other parameters were: $P_0 = 88$ kN, $E = 2.1 \times 10^5$ MPa, $h = 0.003$ m, $R = 0.045$ m, $\rho = 7530$ kg/m$^3$, $\gamma = 2000$ kg/m$^3$, $a = 0.611$ m/ms, $b = 0.357$ m/ms, used later as the reference parametric set. Hereinafter, the parameters of the surrounding medium corresponded to the averaged characteristics of loam.

As seen in Fig. 2, the asymptotic solution (17) accurately describes the finite difference solution in the range between the impact treatment start and the first reflected wave arrival from the external boundary. Reflections from the boundaries $r = R$, $r = R_2$ causes oscillations of (19) about relatively static value.

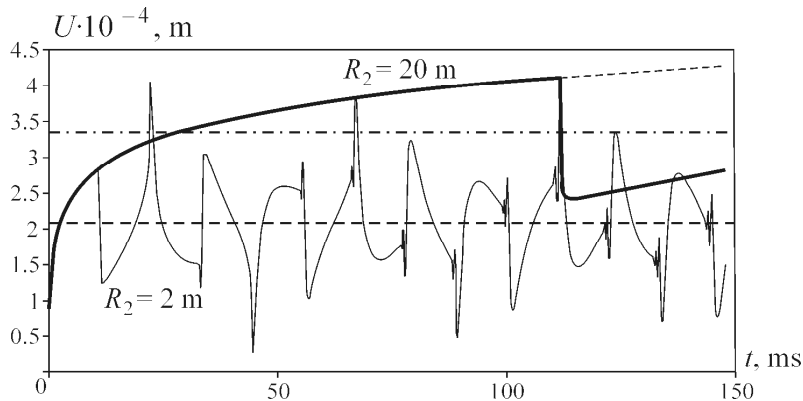

**Fig. 2.** Oscillograms of displacements under step loading.

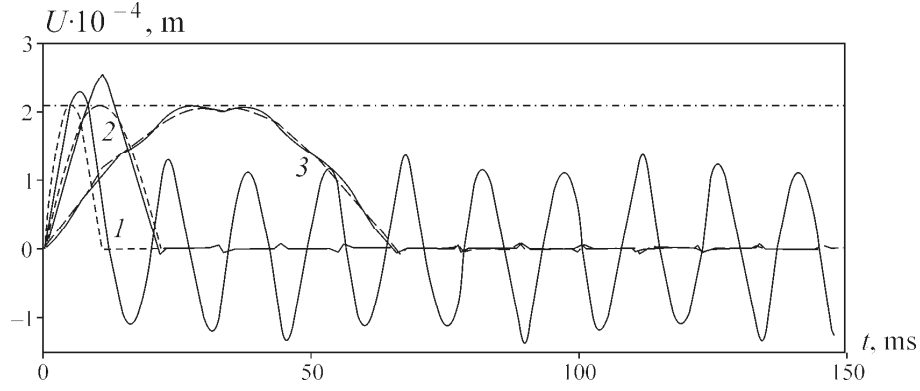

**Fig. 3.** Oscillograms of displacements under pulsed loading.

In case of the pulsed loading (1) of the pipe in the medium with radius $R_2$, the imaged asymptotic solution of (16) has the form of:

$$U^L = \frac{P_0 \omega_*(1+e^{-p\pi/\omega_*})}{S_t\rho(p^2+\omega_*^2)(p^2+\beta^2)}, \quad \beta^2 = \frac{2\pi G}{S_t\rho\ln(R_2/R)}.$$

Converse of the above expression yields time dependence of the pipe displacements:

$$U = \frac{P_0}{S_t\rho\beta(\omega_*^2-\beta^2)}\begin{cases}\omega_*\sin\beta t - \beta\sin\omega_* t, & t\leq t_0,\\ \omega_*(\sin\beta t + \sin\beta(t-t_0)), & t > t_0.\end{cases}$$

If $\omega_* \ll \beta$, then:

$$U = U_0\begin{cases}\sin\omega_* t, & t\leq t_0,\\ 0, & t > t_0.\end{cases} \quad U_0 = \frac{P_0}{S_t\rho\beta^2}. \tag{20}$$

Figure 3 illustrates calculated pulsed loading of a pipe in a medium of fixed $R_2 = 2$ m. the pulse duration is $t_0 = 2(R_2-R)/b$ for curve *1*, $t_0 = 4(R_2-R)/b$ for curve *2* and $t_0 = 12(R_2-R)/b$ for curve *3*; the other input parameters are taken from the reference parametric set. The sold lines show the finite difference solution, the dashed lines—approximated solution (20); dash-and-dotted lines—$U_0$. The comparison of (20) and numerical results yields that if the pulse duration is divisible by time of four runs of the wave between the boundaries of the medium: $t_0 = 4n(R_2-R)/b$ $(n=1,2,...)$, then the waves reflected from the external boundary kill the pulse. As follows from Fig. 3, the pipe displacement is quantitatively and qualitatively correct described with the analytical solution (20).

## 4. PLANE PROBLEM ON ELASTIC INTERACTION AT THE PIPE AND SOIL INTERFACE

The finite difference calculation of the system (2)–(8) was performed for the soil layer with thickness $L = L_1 = 0.2$ m $= 2h_z$ and radii $R_2 = 2$ m and $R_2 = 20$ m under step loading. Analysis of the solutions showed that the logarithmic dependence of the solution on time is typical of the plane problem as with the one-dimensional problem. The qualitative behavior of the both is coincident.

Figure 4 present the numerical calculations of the pile displacements versus time under step loading for different lengths of the pipe and thicknesses of the medium ($L = L_1$). The problem parameters were $R_2 = 20$ m, $h_z = 0.1$ m; the rest—from the reference set. The analysis showed that for a longer pipe, the value of the displacements at the cross-section $z = 0$ at $t = 100$ ms lowers in proportion to the logarithmic function: $F(L) = (898\ln(L)+2061)^{-1}$. For low $L$ (soil layer thickness in the line of $z$ equals two steps of the mesh, $L = 2h_z$), the finite difference solution produces higher error than $F(L)$.

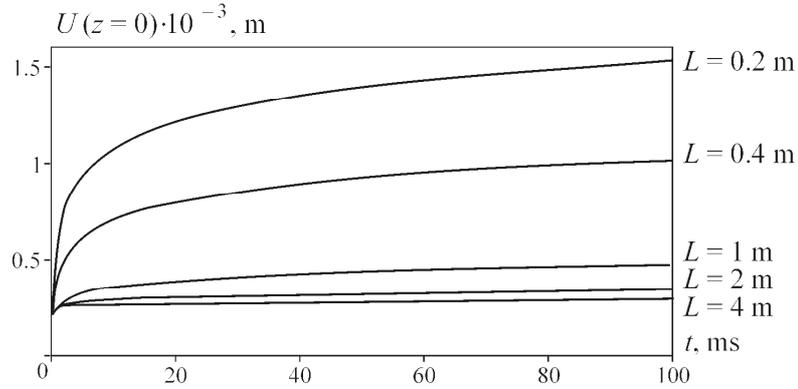

**Fig. 4.** Oscillograms of displacements at various lengths of the pipe under step loading.

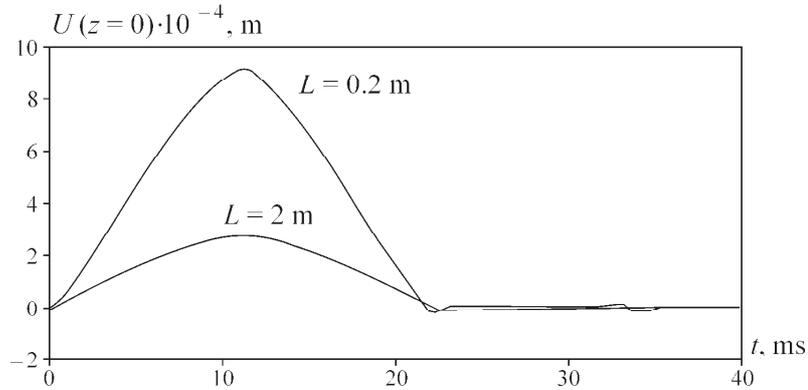

**Fig. 5.** Oscillograms of displacements under pulsed loading.

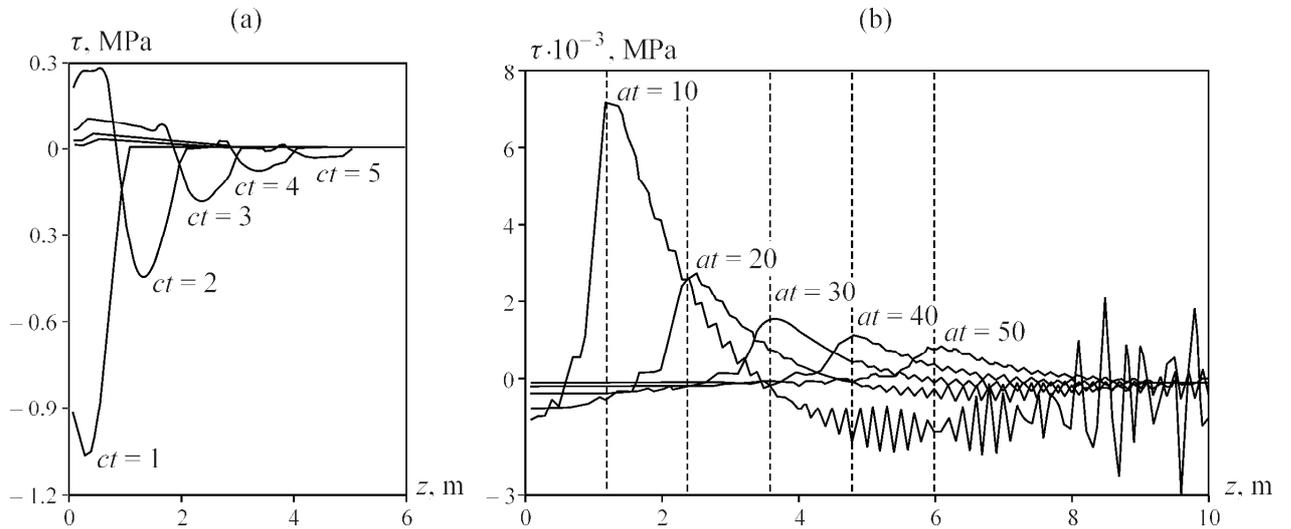

**Fig. 6.** Curves of shear stresses on the pipe surface at different time moments.

The curves in Fig. 5 are shown for the pipe displacements under half-sine pulse (1) with duration $t_0 = 4(R_2 - R)/b$ at two values of $L$ ($L = L_1$). The problem included $R_2 = 2$ m, $h_z = 0.1$ m, and the other parameters taken from the reference set. Qualitatively, the solutions of the plan and one-dimensional problems behave the same way (compare Figs. 3 and 5).

Figure 6 shows curves of the shear stress exerted by the soil on the pipe surface at various time moments. The input parameters were $L = 7.5$ m, $L_1 = 4$ m, $R_2 = 0.8$ m, $t_0 = 0.25$ ms and the rest taken from the reference set. The vertical dashed lines in Fig. 6b show zone of quasi-fronts: $z = at$.

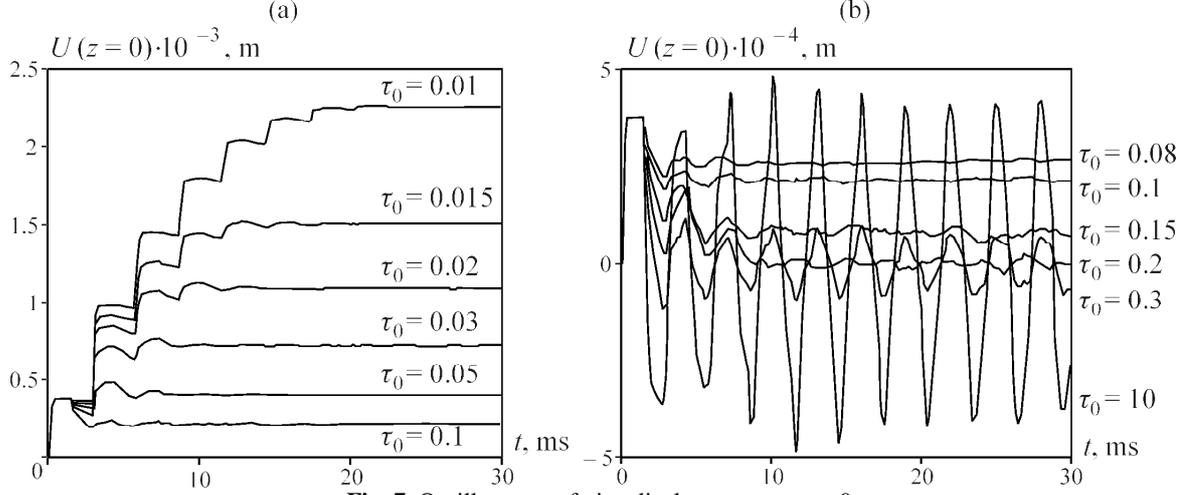

Fig. 7. Oscillograms of pipe displacements at $z = 0$.

In the vicinity of the quasi-front of P-wave with velocity $c$ in the pipe, the maximum shear stress amplitude decreases exponentially as $\sim e^{-0.85ct}$. In the vicinity of the quasi-front of P-wave in the medium the maximum shear stress amplitude decreases as $t^{-4/3}$ (Fig. 6b). In length of time, perturbances that move with the P-wave velocity in soil start dominating. Using Fig. 6b it is possible o evaluate friction force at which the pipe starts slip in soil, depending on the pipe length and the impact pulse amplitude.

## 5. PLANE PROBLEM SOLUTION CONSIDERING SLIP AT PIPE AND SOIL INTERFACE

Figure 7 shows curves of displacement of the pipe end $z = 0$ versus time under half-sine impact pulse and different shear stress values: $\tau_0 = 0.1$, 0.05, 0.03, 0.02, 0.015, 0.01 MPa in Fig. 7a; $\tau_0 = 10$, 0.3, 0.2, 0.15, 0.1, 0.08 MPa in Fig. 7b. The parameters were $L = 7.5$ m, $L_1 = 4$ m, $R_2 = 0.8$ m, $t_0 = 0.22$ ms, and the rest were taken from the reference set. The value $\tau_0 = 10$ MPa conforms with the elastic interaction, the pipe end $z = 0$ oscillates around zero. The amplitude of the oscillations decreases with decrease in the shear stress amplitude (curves $\tau_0 = 0.3$ and $\tau_0 = 0.2$ MPa in Fig. 7b) and residual displacement (curves $\tau_0 = 0.15$ and $\tau_0 = 0.1$ MPa in Fig. 7b and all curves in Fig. 7a) grows as $\tau_0$ becomes lower. As seen in Fig. 7b, when $\tau_0 \leq 2P_0 / LP_t$, the pipe slips in the medium. Analysis of the displacements at $t = 100$ ms versus $\tau_0$ showed that this dependence is close to the inverse function: $\max U(\tau_0) \approx 28\tau_0^{-1}$.

According to the calculated data on displacements of pipes of varied length, the average displacement $2U_{av} = (\max_z U + \min_z U)\big|_{t=100\text{ ms}}$ is in inverse proportion to the pipe length. This result is qualitatively different from the behavior of the pipe elastically embedded in soil (see Fig. 4) for there is the logarithmic dependence in the latter case.

So, displacement of a pipe with dry friction is inverse proportional to friction force $F_f = P_t L \tau_0$.

The calculations illustrated in Fig. 8 were performed for two models of impact pulse loading and different ultimate shear stresses. The thin curve shows model (2)–(10); the thick curve illustrates model (2)–(4), (9) disregarding external medium action from [21]. The problem parameters were $L = 7.5$ m, $L_1 = 4$ m, $R_2 = 0.8$ m, $t_0 = 0.22$ ms, and the rest were taken from the reference set. The two models reach mapping when $\tau_0 < 0.02$ МПа, i.e. when $\tau_0 < P_0/(4LP_t)$. The latter inequality means that when the pipe length is below some definite value governed by the impact pulse amplitude, ultimate shear stress and perimeter of the pipe— $L < L_* = P_0/(4\tau_0 P_t)$ —the exaction of the external medium is possible to neglect and use a simpler model.

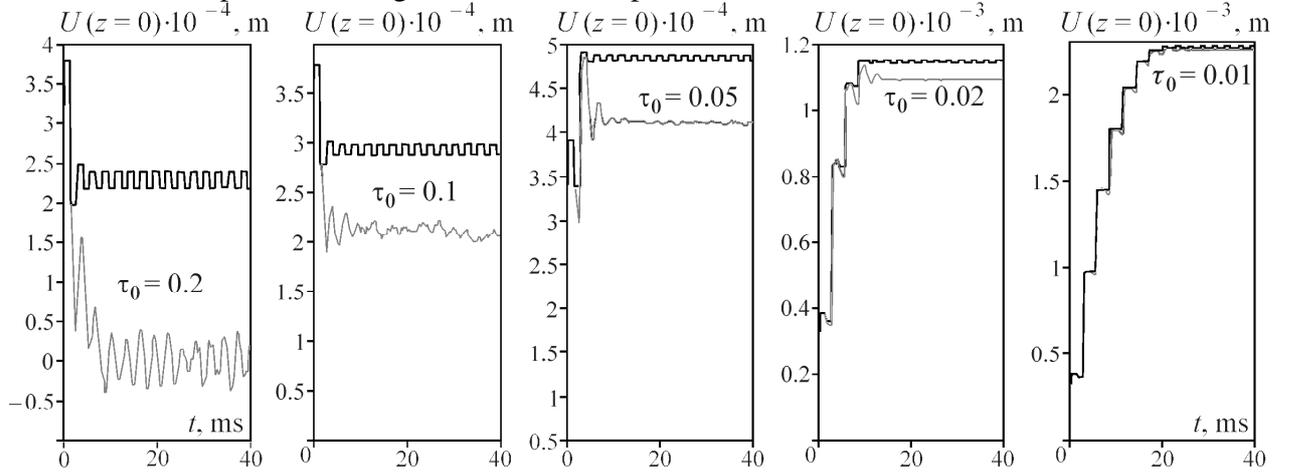

**Fig. 8.** Oscillograms of pipe displacements at $z = 0$.

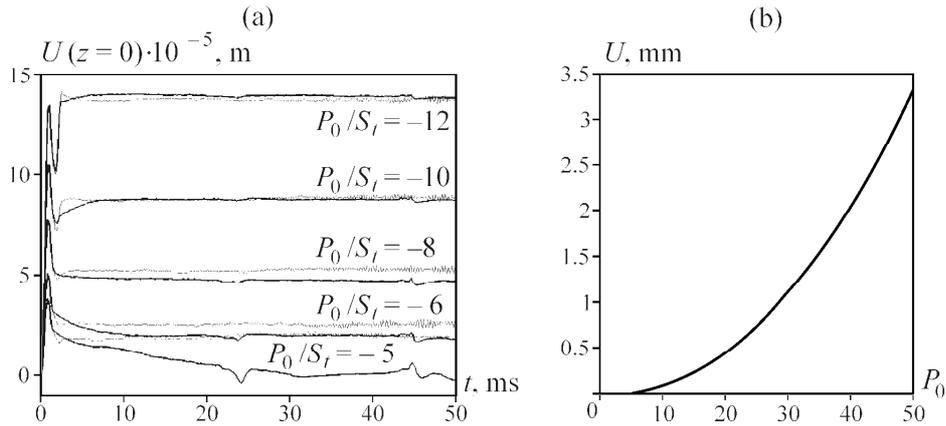

**Fig. 9.** (a) Oscillograms of displacements at different impact pulse amplitude and (b) curve of the displacements and impact pulse amplitude.

The oscillograms of displacements of the pipe end $z = 0$ in Fig. 9a are plotted at varied impact pulse amplitudes. The problem parameters were $L = L_1 = 4$ m, $t_0 = 1$ ms, $E = 2.03 \times 10^5$ MPa, $h = 0.01$ m, $R = 0.1625$ m, $\rho = 7805$ kg/m$^3$, $h_z = 0.1$ m, $\tau_0 = 0.02$ MPa. The thin curves are for $R_2 = 0.2$ m and the thick curves are fore $R_2 = 4$ m. Under small amplitude pulses ($P_0/S_t < 10$ kN) or, which is of the same sort of thin, under low $P_0/F_f$, the external medium radius has influence on the residual displacement. So, when slip begins, the external medium radius has almost no effect of the wave process in the pipe. Figure 9b shows the curves of the pipe displacement at $z = 0$ at $t = 100$ ms and the impact pulse amplitude. The finite difference solution is shown by the thin curve, and its approximated described with the function $U(P_0/S_t) = 0.0014 P_0^2/S_t^2 - 0.0059 P_0/S_t - 0.008$ is displayed by the thick curve. The obtained result qualitatively confirms the deduction on that displacement is proportional to the squared amplitude of the impact pulse.

Figure 10a shows oscillograms of displacements at different pipe lengths ($L = L_1 = 4$, 10, 20, 30 m) and different impact pulse amplitudes, $R_2 = 4$ m, the other parameters are equal to those in Fig. 9. The aim was to find the impact pulse amplitude, such that the average residual displacements are non-zero, i.e. slip start. Figure 10b shows the relationship of the impact pulse amplitude and the pipe length (solid curve) based on the analysis of data in Fig. 10a. The approximation of this relationship is linear: $P_0 / S_t \approx 0.729 L + 1.86$ (dashed curve in Fig. 10b).

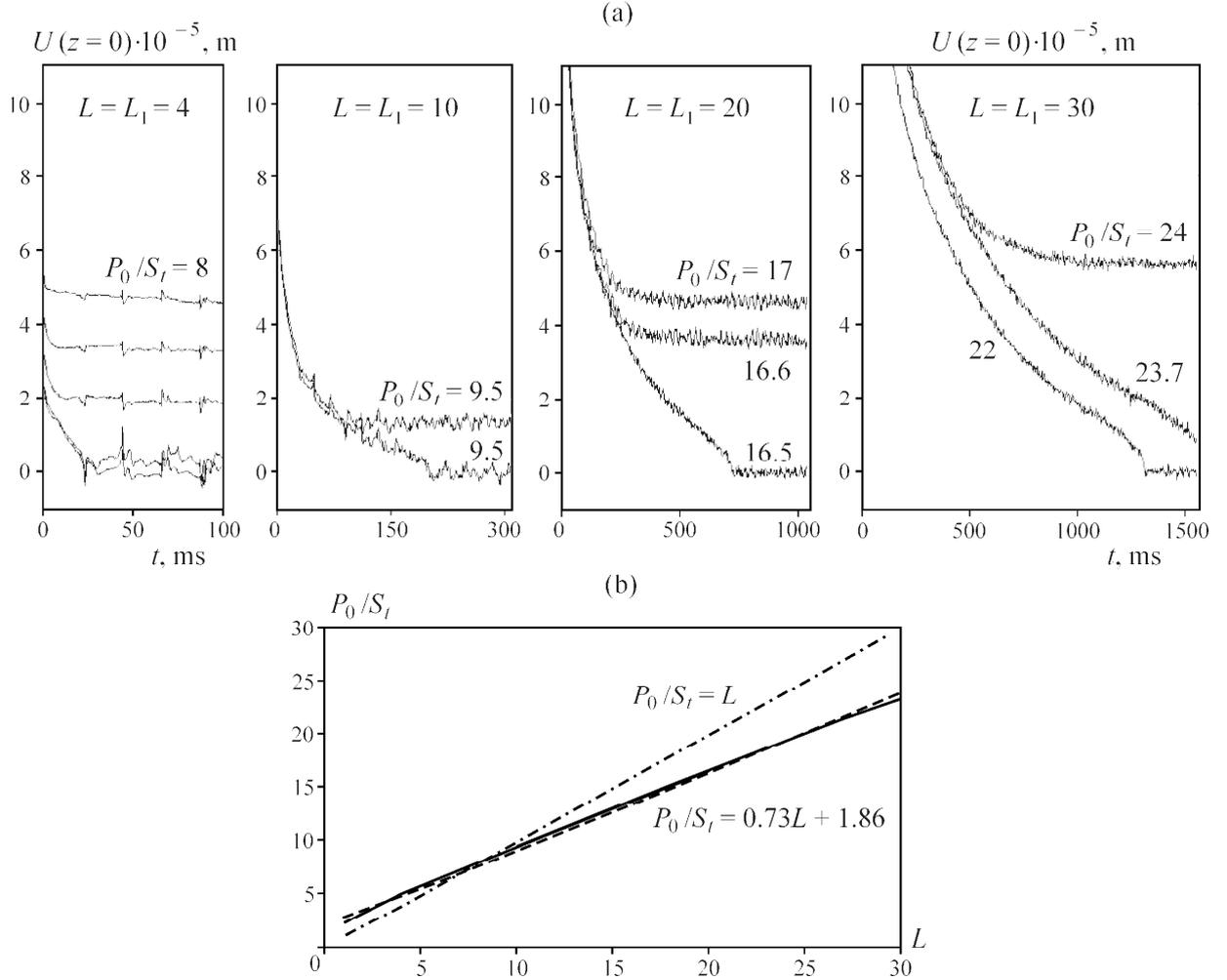

**Fig. 10.** (a) Oscillograms of the pipe displacements under different impact pulse amplitudes; (b) relationship of the slip-initiated impact pulse amplitude and the pipe length.

In [21], using model with neglected formability of soil around a pipe, it was shown that, given the indicated parameters, the impact pulse to reach the pipe ends requires fulfillment of the equality $P_0 / S_t = \tau_0 P_t L / (2 S_t) \approx L$ (dash-dotted curve in Fig. 10b). The calculations with the motion of the external medium taken into account show that when $L > 6$ m, the impact pulse amplitude may be smaller than the above estimate.

## CONCLUSIONS

The author has derived analytical solutions of one-dimensional radial problem on pipe displacements in an elastic host medium under axial impact loading. It is shown that the numerical and analytical solutions of the one-dimensional problem accurately coincide. The analytical solutions qualitatively correctly describe the solution of two-dimensional problem on the elastic interaction of pipe and host medium.

The finite difference calculations of the plain problem on the elastic interaction of the pipe and soil, with dry friction at the pipe and soil interface, performed within a wide range of input parameters, have shown that:

—for the pipe longer than $L_*/4$, where $L_* = P_0/(\tau_0 P_t)$, the soil motion is to be taken into account;

—for the pipe shorter than $L_*/4$, taking into account motion of soil that interacts with the pip[e by the Coulomb law of dry friction weakly influences the results, and the calculation may use simpler models with no regard for the soil motion;

—for the pipe shorter than $2L_*$, slip of the pipe in soli is initiated, i.e. the pipe oscillates elastically in soil;

—the impact pulse amplitude that initates slip of the pipe linearly relates with the pipe length;

—the slip value is in inverse proportion to the friction force $F_f = P_t L \tau_0$.